\documentclass[leqno]{article}

\title{Shape preserving properties of the Bernstein polynomials with integer coefficients}
\author{Borislav R. Draganov}
\date{}

\usepackage{amsmath,amsthm,amsfonts}
\usepackage[mathscr]{eucal}
\usepackage[english]{babel}
\usepackage{fancyhdr}
\usepackage{graphicx}
\usepackage{caption}

\allowdisplaybreaks[3]

\newtheorem{thm}{Theorem}[section]
\newtheorem{prop}[thm]{Proposition}
\newtheorem{cor}[thm]{Corollary}
\newtheorem{lem}[thm]{Lemma}
\theoremstyle{definition}
\newtheorem{defn}[thm]{Definition}
\theoremstyle{remark}
\newtheorem{rem}[thm]{\bf Remark}
\newtheorem{ex}[thm]{\bf Example}

\numberwithin{equation}{section}

\newcommand{\thmref}[1]{Theorem~\ref{#1}}
\newcommand{\propref}[1]{Proposition \ref{#1}}
\newcommand{\corref}[1]{Corollary \ref{#1}}
\newcommand{\lemref}[1]{Lemma \ref{#1}}
\newcommand{\remref}[1]{Remark \ref{#1}}

\newcommand{\N}{\mathbb{N}}
\newcommand{\R}{\mathbb{R}}
\newcommand{\Z}{\mathbb{Z}}
\newcommand{\la}{\left\langle}
\newcommand{\ra}{\right\rangle}

\begin{document}

\maketitle
\bigskip

\thispagestyle{fancy}
\fancyhf{}
\renewcommand{\headrulewidth}{0pt}
\lhead{}
\renewcommand{\footrulewidth}{0.3pt}
\lfoot{\footnotesize This work was supported by grant DN 02/14 of the Fund for Scientific Research of the Bulgarian Ministry of Education and Science.}

\begin{abstract}
The Bernstein polynomials with integer coefficients do not generally preserve monotonicity and convexity. We establish sufficient conditions under which they do. We also observe that they are asymptotically shape preserving.
\end{abstract}

\bigskip
\noindent
{\footnotesize \leftskip25pt \rightskip25pt {\sl  AMS} {\it classification}: 41A10, 41A29, 41A35, 41A36.\\[2pt]
\par}
\bigskip

\section{Main results}

The Bernstein polynomial is defined for $n\in\N_+$, $f\in C[0,1]$ and $x\in [0,1]$ by
\[
B_n f(x):=\sum_{k=0}^n f\left(\frac{k}{n}\right)p_{n,k}(x),\quad p_{n,k}(x):=\binom{n}{k}x^k (1-x)^{n-k}.
\]
It is known that if $f\in C[0,1]$, then (see e.g.~\cite[Chapter 1, Theorem 2.3]{De-Lo:CA})
\[
\lim_{n\to\infty} B_n f(x)=f(x)\quad\text{uniformly on }[0,1].
\]

In order to show that any continuous function on $[0,1]$, which has integer values at the ends of the interval, can be approximated with algebraic polynomials with integer coefficients, Kantorovich \cite{Ka} introduced the operator
\[
\widetilde{B}_n(f)(x):=\sum_{k=0}^n \left[f\left(\frac{k}{n}\right) \binom{n}{k}\right] x^k (1-x)^{n-k},
\]
where $[\alpha]$ denotes the largest integer that is less than or equal to the real $\alpha$. L.~Kantorovich showed that if $f\in C[0,1]$ is such that $f(0),f(1)\in\Z$, then (see \cite{Ka}, or e.g.~, \cite[pp.~3--4]{Fe}, or \cite[Chapter 2, Theorem 4.1]{Lo-Go-Ma:CA})
\[
\lim_{n\to\infty} \widetilde{B}_n(f)(x)=f(x)\quad\text{uniformly on }[0,1].
\]

Instead of the integer part we can take the nearest integer. More precisely, if $\alpha\in\R$ is not a half-integer, we set $\la\alpha\ra$ to be the integer at which the minimum $\min_{m\in\Z} |\alpha-m|$ is attained. When $\alpha$ is a half-integer, we can define $\la\alpha\ra$ to be either of the two neighbouring integers even without following a given rule. The results we will prove are valid regardless of our choice in the latter case. The integer modification of the Bernstein polynomial based on the nearest integer function is given by
\[
\widehat{B}_n(f)(x):=\sum_{k=0}^n \la f\left(\frac{k}{n}\right) \binom{n}{k} \ra x^k (1-x)^{n-k}.
\]
Similarly to \cite{Ka}, it is shown that
\[
\lim_{n\to\infty} \widehat{B}_n(f)(x)=f(x)\quad\text{uniformly on }[0,1]
\]
provided that $f\in C[0,1]$ and $f(0),f(1)\in\Z$.

Let us note that the operators $\widetilde{B}_n$ and $\widehat{B}_n$ are not linear for $n\ge 2$.

As is known, the Bernstein polynomials possess good shape preserving properties. In particular, if $f$ is monotone, then $B_n f$ is monotone of the same type, or, if $f(x)$ is convex or concave so is, respectively, $B_n f$ (see e.g.~\cite[Chapter 10, Theorem 3.3, (i) and (ii)]{De-Lo:CA}). Our main goal is to extend these assertions to the integer forms of the Bernstein polynomials.

The operators $\widetilde{B}_n$ and $\widehat{B}_n$ possess the property of simultaneous approximation, that is, the derivatives of $\widetilde{B}_n(f)$ and $\widehat{B}_n(f)$ approximate the corresponding derivatives of $f$ in the uniform norm on $[0,1]$. This was established in \cite{Dr:SimBernInt,Dr:SimBernIntC} under certain necessary and sufficient conditions, as estimates of the rate the convergence were proved. Hence, trivially, under these conditions, if $f^{(r)}(x)$ is strictly positive or negative, then so are $(\widetilde{B}_n(f))^{(r)}(x)$ and $(\widehat{B}_n(f))^{(r)}(x)$ at least for $n$ large enough, depending on $f$. We will establish sufficient conditions on the shape of $f$ that imply the corresponding monotonicity or convexity of $\widetilde{B}_n(f)$ and $\widehat{B}_n(f)$ for all $n$ regardless of the smoothness of $f$. 

The properties we will present below are not hard to prove. However, they seem interesting and might be useful in the applications of the approximation of functions by polynomials with integer coefficients and in CAGD. Let us note that computer manipulation of polynomials with integer coefficients is faster. 

The operators $\widetilde{B}_n$ and $\widehat{B}_n$ do not generally preserve monotonicity or convexity. We include counter examples in Section 4. It is quite straightforward to show that the monotonicity of $f(x)$ implies the monotonicity of the same type of $\widetilde{B}_n(f)(x)$ and $\widehat{B}_n(f)(x)$ for $n=1$ and $n=2$ (see \remref{n=1,2} below). However, both operators almost preserve monotonicity or convexity. In order to make this precise, we will introduce the notions of asymptotic monotonicity and convexity preservation.

\begin{defn}\label{def1}
Let $X$ be a class of functions defined on $I\subseteq\R$ and $L_n:X\to X$, $n\in\N_+$, be a family of operators. We say that $L_n$ \emph{uniformly asymptotically preserves monotonicity} on $X$ if there exist $n_0\in\N_+$ and functions $\varepsilon_n,\eta_n:I\to\R$, $n\ge n_0$, with the properties:
\begin{enumerate}
\renewcommand{\labelenumi}{(\roman{enumi})}

\item $\lim_{n\to\infty}\varepsilon_n(x)=\lim_{n\to\infty}\eta_n(x)=0$ uniformly on $I$;

\item If $f(x)$ is monotone increasing on $I$, then so is $L_n(f)(x)+\varepsilon_n(x)$ for all $n\ge n_0$;

\item If $f(x)$ is monotone decreasing on $I$, then so is $L_n(f)(x)+\eta_n(x)$ for all $n\ge n_0$.

\end{enumerate}	
	
\end{defn}

\begin{rem}
Let us note that conditions (ii) and (iii) are equivalent if the operators $L_n$ are linear.
\end{rem}

We will show that the following result holds.

\begin{thm}\label{thm2m}
The operators $\widetilde{B}_n$ and $\widehat{B}_n$ uniformly asymptotically preserve monotonicity on the class of continuous functions on $[0,1]$ with integer values at $0$ and $1$.
\end{thm}

Similarly, we introduce the following notion.

\begin{defn}\label{def2}
Let $X$ be a class of functions defined on $I\subseteq\R$ and  $L_n:X\to X$, $n\in\N_+$, be a family of operators. We say that $L_n$ \emph{uniformly asymptotically preserves convexity} on $X$ if there exist $n_0\in\N_+$ and functions $\varepsilon_n,\eta_n:I\to\R$, $n\ge n_0$, with the properties:
	\begin{enumerate}
		\renewcommand{\labelenumi}{(\roman{enumi})}
		
		\item $\lim_{n\to\infty}\varepsilon_n(x)=\lim_{n\to\infty}\eta_n(x)=0$ uniformly on $I$;
		
		\item If $f(x)$ is convex on $I$, then so is $L_n(f)(x)+\varepsilon_n(x)$ for all $n\ge n_0$;

		\item If $f(x)$ is concave on $I$, then so is $L_n(f)(x)+\eta_n(x)$ for all $n\ge n_0$.
			
	\end{enumerate}	
	
\end{defn}

\begin{rem}
As above, if the operators $L_n$ are linear, then (ii) and (iii) are equivalent.
\end{rem}

We will show that $\widetilde{B}_n$ and $\widehat{B}_n$ possess the property described in the definition.

\begin{thm}\label{thm2c}
	The operators $\widetilde{B}_n$ and $\widehat{B}_n$ uniformly asymptotically preserve convexity on the class of continuous functions on $[0,1]$ with integer values at $0$ and $1$.
\end{thm}

On the other hand, it will be useful to establish sufficient conditions on the function $f$ under which we have that $\widetilde{B}_n(f)$ and $\widehat{B}_n(f)$ are monotone, or, respectively, convex or concave. A straightforward corollary of some of our main results is the following assertion.

\begin{thm}\label{thm1m}
	Let $f:[0,1]\to\R$ and $f(0),f(1)\in\Z$.
	\begin{enumerate}
	\renewcommand{\labelenumi}{(\alph{enumi})}
	
		\item	If $f(x)-x$ is monotone increasing on $[0,1]$, then so are $\widetilde{B}_n(f)(x)$ and $\widehat{B}_n(f)(x)$ for all $n$.

		\item If $f(x)+x$ is monotone decreasing on $[0,1]$, then so are $\widetilde{B}_n(f)(x)$ and $\widehat{B}_n(f)(x)$ for all $n$.

\end{enumerate}

\end{thm}

Let us explicitly note that if $f(x)-x$ is monotone increasing on $[0,1]$, then so is $f(x)$, and similarly, if $f(x)+x$ is monotone decreasing on $[0,1]$, then so is $f(x)$.

Also, we will establish the following stronger result.

\begin{thm}\label{thm3m}
	Let $f:[0,1]\to\R$ and $f(0),f(1)\in\Z$. Set for $n\in\N_+$ and $x\in [0,1]$
	\[
\varphi_n(x):=(n+1)\int_0^1 t(1-t)^{n(1-x)}\frac{(1-t)^{nx-1}-t^{nx-1}}{1-2t}\,dt.
	\]
	\begin{enumerate}
	\renewcommand{\labelenumi}{(\alph{enumi})}
	
		\item	If $f(x)-\varphi_{n}(x)$ is monotone increasing on $[0,1]$, then so are $\widetilde{B}_n(f)(x)$ and $\widehat{B}_n(f)(x)$.

		\item If $f(x)+\varphi_{n}(x)$ is monotone decreasing on $[0,1]$, then so are $\widetilde{B}_n(f)(x)$ and $\widehat{B}_n(f)(x)$.

\end{enumerate}

\end{thm}

As it follows from \remref{refvarphi}, the function $\varphi_n(x)$ is monotone increasing on $[0,1]$ for each $n\in\N_+$ and it is of small magnitude---it satisfies the estimates
\[
0\le\varphi_n(x)\le\frac{6}{n},\quad x\in\left[\frac{1}{n},1-\frac{1}{n}\right].
\]

In Section 2 we will establish even less restrictive conditions on $f$ that imply the monotonicity of $\widetilde{B}_n(f)$ and $\widehat{B}_n(f)$. They show how to construct functions $\varphi_n$, which beside the property given in the theorem above, are also such that $|\varphi_n(x)|\le c/n$ for all $x\in [0,1]$ and all $n\in\N_+$, where $c$ as an absolute positive constant; moreover, the functions $\varphi_n$ can be constructed in such a way that if $f(x)\mp\varphi_{n_0}(x)$ is monotone increasing, respectively, decreasing on $[0,1]$ with some $n_0$, then so are $\widetilde{B}_n(f)(x)$ and $\widehat{B}_n(f)(x)$ for all $n\ge n_0$.

Concerning the preservation of convexity and concavity, we will establish

\begin{thm}\label{thm1c}
	Let $f:[0,1]\to\R$ and $f(0),f(1)\in\Z$. Set 
	\[
	\Phi(x):=6\big(x\ln x+(1-x)\ln(1-x)\big).
	\]
	\begin{enumerate}
	\renewcommand{\labelenumi}{(\alph{enumi})}
		\item If $f(x)-\Phi(x)$ is convex on $[0,1]$, then so are $\widetilde{B}_n(f)(x)$ and $\widehat{B}_n(f)(x)$ for all $n$.
		
		\item If $f(x)+\Phi(x)$ is concave on $[0,1]$, then so are $\widetilde{B}_n(f)(x)$ and $\widehat{B}_n(f)(x)$ for all $n$.
	\end{enumerate}

\end{thm}

Note that $\Phi(x)$ is convex and the assumption $f(x)\mp\Phi(x)$ is convex/concave implies that $f(x)$ is convex/concave, respectively.

A less restrictive sufficient condition is given in the following assertion.

\begin{thm}\label{thm3c}
	Let $f:[0,1]\to\R$ and $f(0),f(1)\in\Z$. Set for $n\in\N_+$, $n\ge 3$, and $x\in [0,1]$
\begin{multline*}
	\Phi_n(x):=(n+1)\int_0^1 \left(t^2+(1-t)^2\right)\\
	\times\frac{(nx-3)t^2(1-t)^{n-2} - (nx-2)t^3(1-t)^{n-3}+t^{nx}(1-t)^{n(1-x)}}{(1-2t)^2}\,dt.
\end{multline*}
	\begin{enumerate}
	\renewcommand{\labelenumi}{(\alph{enumi})}
		\item If $f(x)-\Phi_{n}(x)$ is convex on $[0,1]$, then so are $\widetilde{B}_n(f)(x)$ and $\widehat{B}_n(f)(x)$.
		
		\item If $f(x)+\Phi_{n}(x)$ is concave on $[0,1]$, then so are $\widetilde{B}_n(f)(x)$ and $\widehat{B}_n(f)(x)$.
	\end{enumerate}

\end{thm}

As we will establish in \propref{prconv}, the function $\Phi_n(x)$ is convex on $[0,1]$ and it is of small magnitude---it satisfies the estimates
\[
-\frac{4}{n}\le\Phi_n(x)\le\frac{16}{n},\quad x\in \left[\frac{1}{n},1-\frac{1}{n}\right].
\]

In Section 3 we will establish even less restrictive conditions on $f$ that imply the convexity or concavity of $\widetilde{B}_n(f)$ and $\widehat{B}_n(f)$. They show how to construct functions $\Phi_n$, which beside the property given in the theorem above, are also such that $|\Phi_n(x)|\le c/n$ for all $x\in [0,1]$ and all $n\in\N_+$ with some absolute positive constant $c$; moreover, the functions $\Phi_n$ can be constructed in such a way that if $f(x)\mp\Phi_{n_0}(x)$ is convex, respectively, concave on $[0,1]$ with some $n_0$, then so are $\widetilde{B}_n(f)(x)$ and $\widehat{B}_n(f)(x)$ for all $n\ge n_0$.

We proceed to the proof of the results stated above. In the next section we will establish \thmref{thm2m} as well as sufficient conditions that imply the monotonicity of $\widetilde{B}_n(f)(x)$ and $\widehat{B}_n(f)(x)$. In particular, we will get Theorems \ref{thm1m} and \ref{thm3m}. In Section 3 we derive analogues of these results concerning convexity. We present several examples that illustrate the notion of the asymptotic shape preservation and some of the sufficient conditions stated above in Section 4.

\section{Preserving monotonicity}

We set
\begin{align*}
	\tilde b_n(k)&:=\left[f\left(\frac{k}{n}\right)\binom{n}{k} \right]\,\binom{n}{k}^{-1}
	\intertext{and}
	\hat b_n(k)&:=\la f\left(\frac{k}{n}\right)\binom{n}{k} \ra\,\binom{n}{k}^{-1},
\end{align*} 
where $k=0,\dotsc,n$. Then the operators $\widetilde{B}_n$ and $\widehat{B}_n$ can be written respectively in the form
\begin{align*}
	\widetilde{B}_n (f)(x) &=\sum_{k=0}^n \tilde b_n(k)\,p_{n,k}(x)
	\intertext{and}
	\widehat{B}_n (f)(x) &=\sum_{k=0}^n \hat b_n(k)\,p_{n,k}(x).
\end{align*}
For their first derivatives we have (by direct computation, or see \cite{Ma} or \cite[Chapter 10, (2.3)]{De-Lo:CA})
\begin{align}
(\widetilde{B}_n (f))'(x) &=n\sum_{k=0}^{n-1} \left(\tilde b_n(k+1)-\tilde b_n(k)\right)\,p_{n-1,k}(x)\label{eq3}
\intertext{and}
(\widehat{B}_n (f))'(x) &=n\sum_{k=0}^{n-1} \left(\hat b_n(k+1)-\hat b_n(k)\right)\,p_{n-1,k}(x).\label{eq4}
\end{align}

\begin{proof}[Proof of \thmref{thm2m}]
First, let $f(x)$ be monotone increasing on $[0,1]$. We will estimate from below $(\widetilde{B}_n(f))'(x)$.

Since $f(x)$ is increasing on $[0,1/n]$, then
\[
nf\left(\frac{1}{n}\right)\ge nf(0).
\]

We have that $nf(0)\in\Z$; consequently
\[
\left[nf\left(\frac{1}{n}\right)\right]\ge nf(0).
\]
Hence
\begin{equation}\label{eq2.25}
\tilde b_n(1)-\tilde b_n(0) = \frac{1}{n}\left(\left[nf\left(\frac{1}{n}\right)\right]-nf(0)\right)\ge 0.
\end{equation}
Next, using $[\alpha]\le\alpha$, we arrive at
\begin{equation}\label{eq2.26}
\begin{split}
	\tilde b_n(n)-\tilde b_n(n-1) &= f(1)-\frac{1}{n}\left[nf\left(\frac{n-1}{n}\right)\right]\\
	& \ge f(1)-f\left(\frac{n-1}{n}\right)\ge 0.
\end{split}
\end{equation}

Now, let $1\le k\le n-2$, $n\ge 3$. Using the trivial inequalities $\alpha-1\le [\alpha]\le\alpha$, we get
\begin{equation}\label{eq2.27}
\begin{split}
	\tilde b_n(k+1)-\tilde b_n(k)
	&=\left[f\left(\frac{k+1}{n}\right)\binom{n}{k+1} \right]\,\binom{n}{k+1}^{-1}\\
	&\hspace{4cm}-\left[f\left(\frac{k}{n}\right)\binom{n}{k} \right]\,\binom{n}{k}^{-1}\\
	&\ge \left(f\left(\frac{k+1}{n}\right)\binom{n}{k+1} -1 \right)\,\binom{n}{k+1}^{-1}
	- f\left(\frac{k}{n}\right)\\
	&= f\left(\frac{k+1}{n}\right)-f\left(\frac{k}{n}\right)-\binom{n}{k+1}^{-1}.
	\end{split}
\end{equation}
Therefore
\begin{equation}\label{eq2.28}
\tilde b_n(k+1)-\tilde b_n(k)\ge -\binom{n}{k+1}^{-1},\quad 1\le k\le n-2,\ n\ge 3.
\end{equation}

Below we follow the convention that a sum, whose lower index bound is larger than the upper one, is identically $0$. 

We combine \eqref{eq2.25}, \eqref{eq2.26} and \eqref{eq2.28} with \eqref{eq3}, and use the inequality $\binom{n}{k+1}\ge\binom{n}{2}$ for $k=1,\dotsc,n-3$, and the identity $\sum_{k=0}^{n-1} p_{n-1,k}(x)\linebreak[1]\equiv 1$ to arrive at
\begin{align*}
	(\widetilde{B}_n(f))'(x)&\ge n\sum_{k=1}^{n-2} \left(\tilde b_n(k+1)-\tilde b_n(k)\right) p_{n-1,k}(x)\\
	&\ge -\frac{2}{n-1}\sum_{k=1}^{n-3} p_{n-1,k}(x)-(n-1)x^{n-2}(1-x)\\
	&\ge -\frac{2}{n-1}-(n-1)x^{n-2}(1-x).
\end{align*}
We set $\varepsilon_n(x):=2x/(n-1)+x^n/n +x^{n-1}(1-x)$, $n\ge 2$. It satisfies condition (i). Its derivative is $\varepsilon'_n(x)=2/(n-1)+(n-1)x^{n-2}(1-x)$; hence $\varepsilon_n(x)$ satisfies (ii) in Definition \ref{def1} with $n_0=2$.

The case of monotone decreasing functions is reduced to the one of monotone increasing by applying the latter to the function $\bar f(x):=f(1-x)$ and using that $\widetilde{B}_n(f)(x)=\widetilde{B}_n(\bar f)(1-x)$.

The considerations for the operator $\widehat{B}_n$ are quite similar as we use that $|\alpha-\la\alpha\ra|\le 1/2$.
\end{proof}

\begin{rem}\label{n=1,2}
Formula \eqref{eq3} and estimates \eqref{eq2.25} and \eqref{eq2.26} show that if $f(x)$ is monotone increasing on $[0,1]$, then so are $\widetilde{B}_1(f)(x)$ and $\widetilde{B}_2(f)(x)$. Similarly, we have that $\widehat{B}_1(f)(x)$ and $\widehat{B}_2(f)(x)$ are monotone increasing if $f(x)$ is such.
\end{rem}

We proceed to establishing sufficient conditions on $f$ that imply the monotonicity of $\widetilde{B}_n(f)$ and $\widehat{B}_n(f)$. We first consider the operator $\widetilde{B}_n$ and the case of monotone increasing functions.

\begin{prop}\label{thmmtB}
Let $f:[0,1]\to\R$, $f(0),f(1)\in\Z$ and $n\in\N_+$. If $n\ge 3$, let also $\phi_n:[0,1]\to\R$ be such that
\begin{equation}\label{eqmon}
\phi_n\left(\frac{k+1}{n}\right)-\phi_n\left(\frac{k}{n}\right)\ge \binom{n}{k+1}^{-1},\quad k=1,\dotsc,n-2.
\end{equation}
If $f(x)$ is monotone increasing on $[0,1/n]$ and on $[1-1/n,1]$ and if $f(x)-\phi_n(x)$ is monotone increasing on $[1/n,1-1/n]$, then $\widetilde{B}_n(f)(x)$ is monotone increasing on $[0,1]$.
\end{prop}

\begin{proof}
We will show that $\tilde b_n(k+1)-\tilde b_n(k)\ge 0$, $k=0,\dotsc,n-1$. Then, by virtue of $\eqref{eq3}$ we will have 
$(\widetilde{B}_n(f))'(x)\ge 0$ on $[0,1]$.

As we have already established in \eqref{eq2.25} and \eqref{eq2.26}, $\tilde b_n(k+1)-\tilde b_n(k)\ge 0$ for $k=0$ and $k=n-1$.

Let $1\le k\le n-2$, $n\ge 3$. Since $f(x)-\phi_n(x)$ is monotone increasing on $[1/n,1-1/n]$ and $\phi_n(x)$ satisfies \eqref{eqmon}, we have
\begin{equation*}
\begin{split}
f\left(\frac{k+1}{n}\right)-f\left(\frac{k}{n}\right)
&\ge\phi_n\left(\frac{k+1}{n}\right)-\phi_n\left(\frac{k}{n}\right)\\
&\ge \binom{n}{k+1}^{-1}.
\end{split}
\end{equation*}
Then \eqref{eq2.27} implies that $\tilde b_n(k+1)-\tilde b_n(k)\ge 0$, $k=1,\dotsc,n-2$. The proof is completed.
\end{proof}

Clearly, the function $\phi_n(x):=x$ satisfies \eqref{eqmon} for all $n\in\N_+$; hence \thmref{thm1m}(a) follows for the operator $\widetilde{B}_n$. The next corollary contains a less restrictive choice of $\phi_n(x)$. Actually, the function, defined in it, satisfies \eqref{eqmon} as an equality.

\begin{cor}\label{corm}
Let $f:[0,1]\to\R$ and $f(0),f(1)\in\Z$. Let $n\in\N_+$, $n\ge 3$, be fixed. Set
\begin{equation}\label{eqphi}
\phi_n(x):=(n+1)\int_0^1 t(1-t)^{n(1-x)} \frac{(1-t)^{nx} - t^{nx}}{1-2t}\,dt,\quad x\in [0,1].
\end{equation}
If $f(x)-\phi_n(x)$ is monotone increasing on $[0,1]$, then so is $\widetilde{B}_n(f)(x)$.		
\end{cor}

\begin{proof}
The motivation for the definition of $\phi_n(x)$ comes from the following formula, which is derived by the relationship between the beta and gamma functions (see e.g.~\cite{Su} or \cite[Chapter 10, (1.8)]{De-Lo:CA}). We have
\begin{equation}\label{eq2.4}
\begin{split}
	\int_0^1 t^k (1-t)^{n-k} dt
	&=B(k+1,n-k+1)\\
	&=\frac{\Gamma(k+1)\,\Gamma(n-k+1)}{\Gamma(n+2)}\\
	&=\frac{1}{n+1}\,\binom{n}{k}^{-1}.
\end{split}
\end{equation}
Consequently, for $k=1,\dotsc,n-2$ we have
\begin{equation}\label{eq2.1}
\begin{split}
\phi_n\left(\frac{k+1}{n}\right) - \phi_n\left(\frac{k}{n}\right)
&=(n+1)\int_0^1 t^{k+1} (1-t)^{n-(k+1)}dt\\
&=\binom{n}{k+1}^{-1}.
\end{split}
\end{equation}
Thus $\phi_n(x)$ satisfies \eqref{eqmon}.

It remains to observe that $\phi_n(x)$ is differentiable and
\[
\phi'_n(x)=n(n+1)\int_0^1 t^{nx+1}(1-t)^{n(1-x)} \frac{\ln(1-t) - \ln t}{1-2t}\,dt>0,\quad x\in [0,1].
\]	
Therefore $\phi_n(x)$ is monotone increasing on $[0,1]$; hence so is $f(x)$. 

Now, the assertion of the corollary follows from \propref{thmmtB}.	
\end{proof}

\begin{rem}
	The function $\phi_n(x)$, defined in \eqref{eqphi}, can be represented in the following symmetric form
	\[
	\phi_n(x)=\frac{n+1}{2}\int_0^1 \frac{t(1-t)\left((1-t)^{n(1-x)-1} + t^{n(1-x)-1}\right)\left((1-t)^{nx} - t^{nx}\right)}{1-2t}\,dt.
	\]
\end{rem}

Next, we will note the following elementary estimates for the function $\phi_n(x)$, defined in \eqref{eqphi}.

\begin{lem}\label{prm}
	The function $\phi_n(x)$, defined in \eqref{eqphi},  satisfies the estimates
	\[
	0\le\phi_n(x)\le\frac{4}{n},\quad x\in \left[0,1-\frac{1}{n}\right],\quad n\ge 3.
	\] 
\end{lem}

\begin{proof}
As we noted in the proof of \corref{corm}, $\phi_n(x)$ is monotone increasing; hence
\[
\phi_n(0)\le \phi_n(x)
\le \phi_n\left(1-\frac{1}{n}\right),\quad x\in \left[0,1-\frac{1}{n}\right].
\]
Clearly, $\phi_n(0)=0$. 

Next, summing the equalities in \eqref{eq2.1} on $k=1,\dotsc,n-2$, we arrive at
\[
\phi_n\left(1-\frac{1}{n}\right)-\phi_n\left(\frac{1}{n}\right)=\sum_{k=2}^{n-1} \binom{n}{k}^{-1}.
\]
In view of \eqref{eqphi} and \eqref{eq2.4} with $k=1$, we have
\begin{equation*}
\phi_n\left(\frac{1}{n}\right)
=(n+1)\int_0^1 t(1-t)^{n-1} \,dt
=\binom{n}{1}^{-1}.
\end{equation*}

It remains to take into account that $\binom{n}{k}\ge \binom{n}{2}$ for $k=2,\dotsc,n-2$, to deduce that
\begin{align*}
	\phi_n\left(1-\frac{1}{n}\right) 
	& \le \binom{n}{1}^{-1} + (n-3)\binom{n}{2}^{-1} + \binom{n}{n-1}^{-1} \\
	& \le \frac{4}{n}.
\end{align*}
\end{proof}

Rockett \cite[Theorem 1]{Ro} established a neat formula for the sum of the reciprocals of the binomial coefficients.	

Since generally $[-\alpha]\ne -[\alpha]$ and $\la-\alpha\ra\ne -\la\alpha\ra$ (however, $\la\alpha\ra$ is an odd function for some definitions of the nearest integer), the cases of monotone decreasing or concave functions cannot be reduced, respectively, to the cases of increasing or convex functions by considering $-f$ in place of $f$. However, we can swap between increasing and decreasing functions by means of the transformation $\bar f(x):=f(1-x)$. Thus we derive the following sufficient condition concerning the preservation of the monotone decreasing behaviour from \propref{thmmtB}. 

\begin{prop}\label{thmmdtB}
	Let $f:[0,1]\to\R$, $f(0),f(1)\in\Z$ and $n\in\N_+$. If $n\ge 3$, let also $\psi_n:[0,1]\to\R$ be such that
	\begin{equation}\label{eqmon1}
	\psi_n\left(\frac{k+1}{n}\right)-\psi_n\left(\frac{k}{n}\right)\ge\binom{n}{k}^{-1},\quad k=1,\dotsc,n-2.
	\end{equation}
	If $f(x)$ is monotone decreasing on $[0,1/n]$ and on $[1-1/n,1]$ and if $f(x)+\psi_n(x)$ is monotone decreasing on $[1/n,1-1/n]$, then $\widetilde{B}_n(f)(x)$ is monotone decreasing on $[0,1]$.
\end{prop}

The second assertion of \thmref{thm1m} concerning the operator $\widetilde{B}_n$ follows from the last proposition with $\psi_n(x):=x$. A less restrictive $\psi_n$ is defined in the following corollary of \propref{thmmdtB}.

\begin{cor}\label{cormd}
	Let $f:[0,1]\to\R$ and $f(0),f(1)\in\Z$. Let $n\in\N_+$, $n\ge 3$, be fixed. Set
	\begin{equation*}
	\psi_n(x):=(n+1)\int_0^1 t(1-t)^{n(1-x)+1} \frac{(1-t)^{nx-1} - t^{nx-1}}{1-2t}\,dt,\quad t\in [0,1].
	\end{equation*}
	If $f(x)+\psi_n(x)$ is monotone decreasing on $[0,1]$, then so is $\widetilde{B}_n(f)(x)$.		
\end{cor}

\begin{proof}
	The assertion is established similarly to \corref{corm} as instead of \eqref{eq2.1} we show that
	\[
	\psi_n\left(\frac{k+1}{n}\right) - \psi_n\left(\frac{k}{n}\right)
	=\binom{n}{k}^{-1}
	\]
	for $k=1,\dotsc,n-2$. The function $\psi_n(x)$ is monotone increasing.
\end{proof}

\begin{rem}\label{remmd}
	Similarly to \lemref{prm}, it is shown that $\psi_n(x)$, defined in \corref{cormd}, satisfies
	\[
	0\le\psi_n(x)\le \frac{4}{n},\quad x\in \left[\frac{1}{n},1\right].
	\]
\end{rem}

Analogous results hold for the operator $\widehat{B}_n$. They are verified similarly to \propref{thmmtB}, as we use $|\alpha-\la\alpha\ra|\le 1/2$. Let us note that now the assumptions concerning the two types of monotonicity are symmetric unlike those for the operator $\widetilde{B}_n$.

\begin{prop}\label{thmwhB}
	Let $f:[0,1]\to\R$, $f(0),f(1)\in\Z$ and $n\in\N_+$. If $n\ge 3$, let also $\tilde\varphi_n:[0,1]\to\R$ be such that
	\begin{multline}\label{eq2.2}
	\tilde\varphi_n\left(\frac{k+1}{n}\right)-\tilde\varphi_n\left(\frac{k}{n}\right)\ge \frac{1}{2}\left(\binom{n}{k}^{-1}+\binom{n}{k+1}^{-1}\right),\\ k=1,\dotsc,n-2.
	\end{multline}
	
	\begin{enumerate}
		\renewcommand{\labelenumi}{(\alph{enumi})}
	
	\item If $f(x)$ is monotone increasing on $[0,1/n]$ and on $[1-1/n,1]$ and if $f(x)-\tilde\varphi_n(x)$ is monotone increasing on $[1/n,1-1/n]$, then $\widehat{B}_n(f)(x)$ is monotone increasing on $[0,1]$.

	\item If $f(x)$ is monotone decreasing on $[0,1/n]$ and on $[1-1/n,1]$ and if $f(x)+\tilde\varphi_n(x)$ is monotone decreasing on $[1/n,1-1/n]$, then $\widehat{B}_n(f)(x)$ is monotone decreasing on $[0,1]$.
	
	\end{enumerate}

\end{prop}

The assertions of \thmref{thm1m} for $\widehat{B}_n(f)(x)$ follow from the last proposition with $\tilde\varphi_n(x):=x$. 

\begin{rem}\label{refvarphi}
Another function satisfying \eqref{eq2.2} is
\begin{equation}\label{eq2.3}
\tilde\varphi_n(x):=\frac{n+1}{2}\int_0^1 t(1-t)^{n(1-x)}\frac{(1-t)^{nx-1}-t^{nx-1}}{1-2t}\,dt,\quad x\in [0,1].
\end{equation}
As in the previous cases, it is shown that it is differentiable, as
\begin{equation}\label{eq2.29}
\tilde\varphi'_n(x)=\frac{n(n+1)}{2}\int_0^1 t^{nx}(1-t)^{n(1-x)} \frac{\ln(1-t) - \ln t}{1-2t}\,dt,\quad x\in [0,1].
\end{equation}

Consequently, $\tilde\varphi_n(x)$ is monotone increasing on $[0,1]$ and satisfies the estimates
\[
0\le\tilde\varphi_n(x)\le\frac{3}{n} ,\quad x\in\left[\frac{1}{n},1-\frac{1}{n}\right].
\]
\end{rem}

\begin{proof}[Proof of \thmref{thm3m}]
The function $\tilde\varphi_n(x)$, defined in \eqref{eq2.3}, satisfies \eqref{eq2.2}. Then $\varphi_n(x):=2\tilde\varphi_n(x)$ satisfies the conditions \eqref{eqmon}, \eqref{eqmon1} and \eqref{eq2.2}. Since $f(x)-\varphi_{n}(x)$ is monotone increasing on $[0,1]$, then so is $f(x)$. Now, Propositions \ref{thmmtB} and \ref{thmwhB}(a) yield that $\widetilde{B}_n(f)(x)$ and $\widehat{B}_n(f)(x)$ are monotone increasing on $[0,1]$. The proof of assertion (b) of the theorem is similar.
\end{proof}

\section{Preserving convexity}

For the second derivatives of $\widetilde{B}_n (f)$ and $\widehat{B}_n (f)$ we have (by direct computation, or see \cite{Ma} or \cite[Chapter 10, (2.3)]{De-Lo:CA})
\begin{align}
(\widetilde{B}_n (f))''(x) &=n(n-1)\sum_{k=0}^{n-2} \left(\tilde b_n(k+2)-2\tilde b_n(k+1)+\tilde b_n(k)\right)\,p_{n-2,k}(x)\label{eq5}
\intertext{and}
(\widehat{B}_n (f))''(x) &=n(n-1)\sum_{k=0}^{n-2} \left(\hat b_n(k+2)-2\hat b_n(k+1)+\hat b_n(k)\right)\,p_{n-2,k}(x).\label{eq6}
\end{align}

\begin{proof}[Proof of \thmref{thm2c}]
Similarly to the proof of the corresponding result in the monotone case, we estimate the second derivative of $\widetilde{B}_n(f)(x)$ and $\widehat{B}_n(f)(x)$. We will consider in detail only the former operator in the case of convex functions; the arguments for the latter operator are quite alike. The case of concave functions is analogous too.

Let $f(x)$ be convex on the interval $[0,1]$. Then
\[
f\left(\frac{k+2}{n}\right)-2f\left(\frac{k+1}{n}\right)+f\left(\frac{k}{n}\right)\ge 0,\quad k=0,\dotsc,n-2,\ n\ge 2.
\]
Using $\alpha-1\le [\alpha]\le\alpha$ and $f(0)\in\Z$, we get
\begin{equation}\label{eq3.51}
\begin{split}
	&\tilde b_n(2)-2\tilde b_n(1)+\tilde b_n(0)\\
	&\quad =\left[f\left(\frac{2}{n}\right)\binom{n}{2} \right]\,\binom{n}{2}^{-1}
	-2\left[f\left(\frac{1}{n}\right)\binom{n}{1} \right]\,\binom{n}{1}^{-1}+f(0)\\
	&\quad \ge\left(f\left(\frac{2}{n}\right)\binom{n}{2} -1 \right)\,\binom{n}{2}^{-1}
-2f\left(\frac{1}{n}\right)+f(0)\\ 
	&\quad =f\left(\frac{2}{n}\right)-2f\left(\frac{1}{n}\right)+f(0)-\binom{n}{2}^{-1}.
\end{split}	
\end{equation}
Similarly, we get
\begin{multline}\label{eq3.52}
\tilde b_n(1)-2\tilde b_n(n-1)+\tilde b_n(n-2)\\
\ge f(1)-2f\left(\frac{n-1}{n}\right)+f\left(\frac{n-2}{n}\right)-\binom{n}{n-2}^{-1},\quad n\ge 2.
\end{multline}

Let $k=1,\dotsc,n-3$, $n\ge 4$. Just analogously, we arrive at the estimates
\begin{equation}\label{eq3.53}
\begin{split}
	&\tilde b_n(k+2)-2\tilde b_n(k+1)+\tilde b_n(k)\\
	&\quad =\left[f\left(\frac{k+2}{n}\right)\binom{n}{k+2} \right]\,\binom{n}{k+2}^{-1}\\
	&\qquad\qquad\qquad\qquad -2\left[f\left(\frac{k+1}{n}\right)\binom{n}{k+1} \right]\,\binom{n}{k+1}^{-1}\\
	&\qquad\qquad\qquad\qquad +\left[f\left(\frac{k}{n}\right)\binom{n}{k} \right]\,\binom{n}{k}^{-1}\\
	&\quad \ge\left(f\left(\frac{k+2}{n}\right)\binom{n}{k+2} -1 \right)\,\binom{n}{k+2}^{-1} -2f\left(\frac{k+1}{n}\right)\\
	&\qquad\qquad\qquad\qquad +\left(f\left(\frac{k}{n}\right)\binom{n}{k} -1\right)\,\binom{n}{k}^{-1}\\ 
	&\quad \ge f\left(\frac{k+2}{n}\right) -2f\left(\frac{k+1}{n}\right) + f\left(\frac{k}{n}\right)\\
	&\qquad\qquad\qquad\qquad -\left(\binom{n}{k}^{-1} + \binom{n}{k+2}^{-1}\right).
\end{split}
\end{equation}

Thus we have shown that
\begin{equation*}
\begin{split}
&\tilde b_n(2)-2\tilde b_n(1)+\tilde b_n(0)\ge -\binom{n}{2}^{-1},\\
&\tilde b_n(k+2)-2\tilde b_n(k+1)+\tilde b_n(k)\ge -\left(\binom{n}{k}^{-1} + \binom{n}{k+2}^{-1}\right),\\
&\hspace{7cm} k=1,\dotsc,n-3,\\
&\tilde b_n(1)-2\tilde b_n(n-1)+\tilde b_n(n-2)\ge -\binom{n}{n-2}^{-1}.
\end{split}
\end{equation*}
Consequently, by virtue of \eqref{eq5}, the inequality $\binom{n}{k}\ge\binom{n}{3}$ for $k=3,\dotsc,n-3$, and the identity $\sum_{k=0}^{n-2} p_{n-2,k}(x)\equiv 1$, we have for $n\ge 6$
\begin{align*}
(\widetilde{B}_n (f))''(x) &\ge - n(n-1)\binom{n}{2}^{-1}(1-x)^{n-2}\\
&\qquad -n(n-1)\sum_{k=1}^{n-3}\left(\binom{n}{k}^{-1} + \binom{n}{k+2}^{-1}\right)p_{n-2,k}(x)\\
&\qquad -n(n-1)\binom{n}{n-2}^{-1}x^{n-2}\\
&\ge -2(1-x)^{n-2} - (n-1)(n-2)x(1-x)^{n-3}\\
&\qquad - (n-2)(n-3)x^2(1-x)^{n-4} - n(n-1)\binom{n}{3}^{-1}\sum_{k=3}^{n-3} p_{n-2,k}(x)\\
&\qquad - n(n-1)\binom{n}{3}^{-1}\sum_{k=1}^{n-5} p_{n-2,k}(x) - (n-2)(n-3)x^{n-4}(1-x)^2\\
&\qquad - (n-1)(n-2)x^{n-3}(1-x) - 2x^{n-2}\\
&\ge -2(1-x)^{n-2} - (n-1)(n-2)x(1-x)^{n-3}\\
&\qquad - (n-2)(n-3)x^2(1-x)^{n-4} - \frac{12}{n-2}\\
&\qquad - (n-2)(n-3)x^{n-4}(1-x)^2 - (n-1)(n-2)x^{n-3}(1-x)\\
&\qquad - 2x^{n-2}=:-b_n(x).
\end{align*}
We set
\begin{multline*}
\varepsilon_n(x) := \frac{6x^2}{n-2}-\frac{2(n-3)}{n(n-1)}\big(x^n+(1-x)^n\big)\\
 +\frac{4}{n-1}\left(x^{n-1}+(1-x)^{n-1}\right) +x^{n-2}(1-x)+x(1-x)^{n-2}.
\end{multline*}
We have $\varepsilon_n''(x)=b_n(x)$; hence $\varepsilon_n(x)$ satisfies condition (ii) in Definition \ref{def2} with $n_0=6$. Clearly, it satisfies condition (i) too.
\end{proof}

Further, we will derive sufficient conditions on the function $f$ that imply the convexity and concavity of $\widetilde{B}_n(f)(x)$ and $\widehat{B}_n(f)(x)$.

\begin{prop}\label{thmctB}
Let $f:[0,1]\to\R$ and $f(0),f(1)\in\Z$. Let $n\in\N_+$, $n\ge 2$, be fixed and $\Phi_n:[0,1]\to\R$ be such that
\begin{align*}
&\Phi_n\left(\frac{2}{n}\right)-2\Phi_n\left(\frac{1}{n}\right)+\Phi_n(0)
\ge \binom{n}{2}^{-1}\\
&\Phi_n\left(\frac{k+2}{n}\right)-2\Phi_n\left(\frac{k+1}{n}\right)+\Phi_n\left(\frac{k}{n}\right)\\
&\hspace{2cm}\ge \binom{n}{k}^{-1}+\binom{n}{k+2}^{-1},\quad k=1,\dotsc,n-3,\ n\ge 4,\notag
\intertext{and}
&\Phi_n(1)-2\Phi_n\left(\frac{n-1}{n}\right)+\Phi_n\left(\frac{n-2}{n}\right)
\ge \binom{n}{n-2}^{-1}.
\end{align*}
If $f(x)-\Phi_n(x)$ is convex on $[0,1]$, then so is $\widetilde{B}_n(f)(x)$.
\end{prop}

\begin{proof}
Since $f(x)-\Phi_n(x)$ is convex on $[0,1]$, then
\begin{multline*}
	f\left(\frac{k+2}{n}\right)-2f\left(\frac{k+1}{n}\right)+f\left(\frac{k}{n}\right)\\
	\ge \Phi_n\left(\frac{k+2}{n}\right)-2\Phi_n\left(\frac{k+1}{n}\right)+\Phi_n\left(\frac{k}{n}\right),\quad k=0,\dotsc,n-2.
\end{multline*}
Then \eqref{eq3.51}-\eqref{eq3.53} and the assumptions on $\Phi_n(x)$ imply
\[
\tilde b_n(k+2)-2\tilde b_n(k+1)+\tilde b_n(k)\ge 0,\quad k=0,\dotsc,n-2,
\]
which, by virtue of \eqref{eq5}, completes the proof of the proposition.
\end{proof}

Similarly to \propref{thmctB} we prove the following sufficient condition for preserving concavity.

\begin{prop}\label{thmcvtB}
Let $f:[0,1]\to\R$ and $f(0),f(1)\in\Z$. Let $n\in\N_+$, $n\ge 2$, be fixed and $\Phi_n:[0,1]\to\R$ be such that
\begin{equation*}
\Phi_n\left(\frac{k+2}{n}\right)-2\Phi_n\left(\frac{k+1}{n}\right)+\Phi_n\left(\frac{k}{n}\right)\ge 2\binom{n}{k+1}^{-1},\quad k=0,\dotsc,n-2.
\end{equation*}
If $f(x)+\Phi_n(x)$ is concave on $[0,1]$, then so is $\widetilde{B}_n(f)(x)$.
\end{prop}

Similarly to Propositions \ref{thmctB} and \ref{thmcvtB} , we have the following result for the other integer modification of the Bernstein polynomials, the operator $\widehat{B}_n$.

\begin{prop}\label{thmchB}
Let $f:[0,1]\to\R$ and $f(0),f(1)\in\Z$. Let $n\in\N_+$, $n\ge 2$, be fixed and $\Phi_n:[0,1]\to\R$ be such that
\begin{align*}
&\Phi_n\left(\frac{2}{n}\right)-2\Phi_n\left(\frac{1}{n}\right)+\Phi_n(0)
\ge \frac{1}{2}\left(2\binom{n}{1}^{-1}+\binom{n}{2}^{-1}\right),\\
&\Phi_n\left(\frac{k+2}{n}\right)-2\Phi_n\left(\frac{k+1}{n}\right)+\Phi_n\left(\frac{k}{n}\right)\\
&\qquad\ge \frac{1}{2}\left(\binom{n}{k}^{-1}+2\binom{n}{k+1}^{-1}+\binom{n}{k+2}^{-1}\right),\quad k=1,\dotsc,n-3,\ n\ge 4,
\intertext{and}
&\Phi_n(1)-2\Phi_n\left(\frac{n-1}{n}\right)+\Phi_n\left(\frac{n-2}{n}\right)
\ge \frac{1}{2}\left(\binom{n}{n-2}^{-1}+2\binom{n}{n-1}^{-1}\right).
\end{align*}

\begin{enumerate}
\renewcommand{\labelenumi}{(\alph{enumi})}
	
\item If $f(x)-\Phi_n(x)$ is convex on $[0,1]$, then so is $\widehat{B}_n(f)(x)$.

\item If $f(x)+\Phi_n(x)$ is concave on $[0,1]$, then so is $\widehat{B}_n(f)(x)$.

\end{enumerate}
\end{prop}

We proceed to the proof of \thmref{thm1c}.

\begin{proof}[Proof of \thmref{thm1c}]
For $n=1$ the assertion is trivial since $\widetilde{B}_1(f)(x)$ and $\widehat{B}_1(f)(x)$ are linear functions. Let $n\ge 2$. We will verify that the function $\Phi(x)$ defined in the theorem satisfies the conditions in the propositions stated so far in this section. We set
\[
\Delta(k):=\Phi\left(\frac{k+2}{n}\right)-2\Phi\left(\frac{k+1}{n}\right)+\Phi\left(\frac{k}{n}\right),\quad k=0,\dotsc,n-2.
\]

First, we observe that
\begin{align}
	&2\binom{n}{1}^{-1}\ge \frac{1}{2}\left(2\binom{n}{1}^{-1}+\binom{n}{2}^{-1}\right)\ge\binom{n}{2}^{-1},\label{eq3.6}\\
	&2\binom{n}{n-1}^{-1}\ge  \frac{1}{2}\left(\binom{n}{n-2}^{-1}+2\binom{n}{n-1}^{-1}\right)\ge\binom{n}{n-2}^{-1},\label{eq3.7}\\
	&\binom{n}{k}^{-1}+\binom{n}{k+2}^{-1}
	\ge \frac{1}{2}\left(\binom{n}{k}^{-1}+2\binom{n}{k+1}^{-1}+\binom{n}{k+2}^{-1}\right)\label{eq3.8}\\
	&\hspace{7cm} k=1,\dotsc,n-3,\ n\ge 4,\notag
	\intertext{and}
	&\binom{n}{k}^{-1}+\binom{n}{k+2}^{-1}	\ge 2\binom{n}{k+1}^{-1},\quad k=1,\dotsc,n-3,\ n\ge 4.\label{eq3.9}
\end{align}
Relations \eqref{eq3.6} and \eqref{eq3.7} are identical and trivial. It is straightforward to see that \eqref{eq3.8} and \eqref{eq3.9} are equivalent too. Let us verify the last one. It reduces to
\[
(n-k-1)(n-k) + (k+1)(k+2)\ge 2(k+1)(n-k-1).
\]
We divide both sides of the inequality above by $(k+1)(n-k-1)$, to arrive at
\[
\frac{n-k}{k+1} + \frac{k+2}{n-k-1}\ge 2.
\]
It remains to observe that the second term on the left hand-side is larger than the reciprocal of the first one and then to take into account that the sum of a positive real and its reciprocal is always at least $2$.

Thus to show that $\Phi(x)$ satisfies the assumptions in Propositions \ref{thmctB}, \ref{thmcvtB} and \ref{thmchB}, it is sufficient to prove that
\begin{align}
	\Delta(k)&\ge 2\binom{n}{1}^{-1}=\frac{2}{n},\quad k=0,n-2,\label{eq3.4}
	\intertext{and}
	\Delta(k)&\ge \binom{n}{k}^{-1}+\binom{n}{k+2}^{-1},\quad k=1,\dotsc,n-3,\ n\ge 4.\label{eq3.5}
\end{align}

The function $\Phi(x)$ is twice continuously differentiable in $(0,1)$ and
\[
\Phi''(x)=\frac{6}{x(1-x)}.
\]
By Taylor's formula we get for $k=0,\dotsc,n-2$
\begin{equation}\label{eq3.T}
\Delta(k)=\int_{k/n}^{(k+2)/n} M_{n,k}(t)\,\Phi''(t)\,dt,
\end{equation}
where
\[
M_{n,k}(t):=\begin{cases}
t-\frac{k}{n}, &t\in \left[\frac{k}{n},\frac{k+1}{n}\right],\\
\frac{k+2}{n}-t, &t\in \left(\frac{k+1}{n},\frac{k+2}{n}\right].
\end{cases}
\]
For $k=0$ formula \eqref{eq3.T} implies
\begin{align*}
	\Delta(0)&=6\int_0^{1/n} \frac{dt}{1-t} + 6\int_{1/n}^{2/n} \left(\frac{2}{n}-t\right)\frac{dt}{t(1-t)}\\
	&\ge 6\int_0^{1/n} \frac{dt}{1-t}\\
	&\ge\frac{6}{n}.
\end{align*}
Thus \eqref{eq3.4} is verified for $k=0$. The case $k=n-2$ is symmetric to $k=0$.

For  $k=1,\dotsc,n-3$, \eqref{eq3.T} yields
\begin{align*}
	\Delta(k)&\ge\frac{6}{\max_{x\in  [k/n,(k+2)/n]}x(1-x)} \int_{k/n}^{(k+2)/n} M_{n,k}(t)\,dt\\
	&=\frac{6}{n^2\max_{x\in  [k/n,(k+2)/n]}x(1-x)}.
\end{align*}
If $(k+2)/n\le 1/2$, then $\max_{x\in  [k/n,(k+2)/n]}x(1-x)=(k+2)(n-k-2)/n^2$ and $\binom{n}{k}\le\binom{n}{k+2}$; hence \eqref{eq3.5} will follow from
\[
\frac{6}{(k+2)(n-k-2)}\ge 2\binom{n}{k}^{-1}.
\]
This inequality follows from
\[
\frac{3}{n(k+2)}\ge \binom{n}{k}^{-1},
\]
which is trivial for $k=1$, and otherwise follows from
\[
\frac{3}{n(k+2)}\ge \binom{n}{2}^{-1}=\frac{2}{n(n-1)}.
\]

The case $k/n\ge 1/2$ is symmetric to the case just considered.

It remains to verify \eqref{eq3.5} for $k$ such that $1/2\in (k/n,(k+2)/n)$. Then $\max_{x\in  [k/n,(k+2)/n]}x(1-x)=1/4$. The condition $1/2\in (k/n,(k+2)/n)$ is equivalent to $n/2-2<k<n/2$.

If $n$ is even, then $k=n/2-1$ and $\binom{n}{k}=\binom{n}{k+2}$. In this case \eqref{eq3.5} will follow from
\[
\frac{12}{n^2}\ge \binom{n}{n/2-1}^{-1}.
\]
This is verified directly for $n=4$; otherwise, it follows from
\begin{equation}\label{eq3.10}
\frac{12}{n^2}\ge \binom{n}{2}^{-1}=\frac{2}{n(n-1)},
\end{equation}
which is trivial.

Finally, if $n$ is odd, then $k=(n-3)/2$ or $k=(n-1)/2$. These two cases are symmetric and it suffices to consider $k=(n-3)/2$. Then $\binom{n}{k}<\binom{n}{k+2}$. Therefore \eqref{eq3.5} will follow from
\[
\frac{12}{n^2}\ge \binom{n}{(n-3)/2}^{-1}.
\]
This is checked directly for $n=5$; otherwise, it follows from \eqref{eq3.10}.
\end{proof}

We proceed to the proof of \thmref{thm3c}.

\begin{proof}[Proof of \thmref{thm3c}]
Direct computations and \eqref{eq2.4} yield for $n\ge 3$ and $k=0,\dotsc,n-2$ the relation
\begin{equation}\label{eq3.41}
	\Phi_n\left(\frac{k+2}{n}\right)-2\Phi_n\left(\frac{k+1}{n}\right)+\Phi_n\left(\frac{k}{n}\right)
	=\binom{n}{k}^{-1} + \binom{n}{k+2}^{-1}.
\end{equation}
Therefore, by virtue of \eqref{eq3.8} and \eqref{eq3.9}, the function $\Phi_n(x)$ satisfies the conditions in Propositions \ref{thmctB}-\ref{thmchB}; hence the assertions of the theorem follow.
\end{proof}

\begin{prop}\label{prconv}
The function $\Phi_n(x)$, defined in \thmref{thm3c}, is convex on $[0,1]$ and satisfies the estimates

\begin{equation}\label{eq3.43}
-\frac{4}{n}\le\Phi_n(x)\le\frac{16}{n} ,\quad x\in \left[\frac{1}{n},1-\frac{1}{n}\right].
\end{equation}
\end{prop}

\begin{proof}
As we assumed in the statement of \thmref{thm3c}, $n\ge 3$. The function $\Phi_n(x)$ is twice continuously differentiable on $[0,1]$, as
\begin{align*}
	\Phi'_n(x)&=n(n+1)\int_0^1 \left(t^2+(1-t)^2\right)\\
	&\qquad\times\frac{t^2(1-t)^{n-2} - t^3(1-t)^{n-3}+t^{nx}(1-t)^{n(1-x)}(\ln t-\ln(1-t))}{(1-2t)^2}\,dt
	\intertext{and}
	\Phi''_n(x)&=n^2(n+1)\int_0^1 \left(t^2+(1-t)^2\right) t^{nx}(1-t)^{n(1-x)}\left(\frac{\ln t-\ln(1-t)}{1-2t}\right)^2 dt.
\end{align*}
We have that $\Phi''_n(x)>0$ on $[0,1]$; hence $\Phi_n(x)$ is convex on $[0,1]$.

Further, since $\Phi_n(x)$ is convex, then
\begin{equation}\label{eq3.44}
\Phi_n(x)\le\max\left\{\Phi_n\left(\frac{1}{n}\right),\Phi_n\left(1-\frac{1}{n}\right)\right\},\quad x\in \left[\frac{1}{n},1-\frac{1}{n}\right].
\end{equation}
Straightforward computations and \eqref{eq2.4} show that
\begin{equation}\label{eq3.45}
\Phi_n\left(\frac{1}{n}\right)=\binom{n}{1}^{-1} + \binom{n}{3}^{-1}\le\frac{4}{n}.
\end{equation}

The definition of $\Phi_n(x)$ readily yields that $\Phi_n(2/n)=\Phi_n(3/n)=0$. This, combined with \eqref{eq3.44} and \eqref{eq3.45}, implies \eqref{eq3.43} for $n=3,4$. 

To estimate $\Phi_n(1-1/n)$ for $n\ge 5$ we sum relations \eqref{eq3.41} on $k=2,\dotsc,j$ and then on $j=2,\dotsc,n-3$. As we take into account $\Phi_n(2/n)=\Phi_n(3/n)=0$, we arrive at
\begin{equation}\label{eq3.42}
	\Phi_n\left(\frac{n-1}{n}\right)
	=\sum_{k=2}^{n-3} (n-k-2)\binom{n}{k}^{-1} + \sum_{k=2}^{n-3} (n-k-2)\binom{n}{k+2}^{-1}.
\end{equation}

Next, we estimate the right-hand-side of \eqref{eq3.42}:
\begin{align*}
	\sum_{k=2}^{n-3} (n-k-2)\binom{n}{k}^{-1}
	&\le \sum_{k=2}^{n-2} (n-k+1)\frac{k!\,(n-k)!}{n!}\\
	&=(n+1)\sum_{k=2}^{n-2} \binom{n+1}{k}^{-1}\\
	&\le (n+1)\binom{n+1}{2}^{-1} + (n+1)(n-4)\binom{n+1}{3}^{-1}\\
	&\le\frac{8}{n}.
\end{align*}
Similarly, we get
\[
\sum_{k=2}^{n-3} (n-k-2)\binom{n}{k+2}^{-1}\le\frac{8}{n}.
\]
By virtue of the last two estimates, the fact that $\Phi_n(2/n)=\Phi_n(3/n)=0$ and \eqref{eq3.42}, we arrive at
\[
\Phi_n\left(\frac{n-1}{n}\right)\le\frac{16}{n}.
\]
This along with \eqref{eq3.44} and \eqref{eq3.45} imply the upper estimate in \eqref{eq3.43} for $n\ge 5$.

In order to verify the lower estimate, we use that $\Phi_n(x)$ is convex and $\Phi_n(2/n)=\Phi_n(3/n)=0$ to deduce that $\Phi_n(x)$ attains its global minimum on the interval $(2/n,3/n)$. Since $\Phi_n(x)$ is convex, its graph on the interval $[2/n,3/n]$ lies above the secant line through the points $(1/n,\Phi_n(1/n))$ and $(2/n,\Phi_n(2/n))$. Thus we arrive at
\[
\Phi_n(x)\ge \Phi_n\left(\frac{1}{n}\right)(2-nx)\ge -\Phi_n\left(\frac{1}{n}\right),\quad x\in \left[\frac{2}{n},\frac{3}{n}\right].
\]
Hence, taking into account \eqref{eq3.45}, we get the left inequality in \eqref{eq3.43}.
\end{proof}

\section{Exampes}

We will give several examples to illustrate some of the results obtained above.
 
We begin with an example, which shows that the operator $\widetilde{B}_n$ does not preserve monotonicity for all $n$. It can be shown that if $f$ is monotone increasing, then so is $\widetilde{B}_n(f)$ for $n\le 5$. Here is a counterexample for $n=6$. 

\begin{ex}
Let
\[
\begin{aligned}
	&f(0)=0;\quad &f\left(\frac{1}{6}\right)=\frac{50}{60};\quad &f\left(\frac{2}{6}\right)=\frac{56}{60};
	\quad &f\left(\frac{3}{6}\right)=\frac{57}{60};\\
	&f\left(\frac{4}{6}\right)=\frac{58}{60};\quad &f\left(\frac{5}{6}\right)=\frac{59}{60};\quad &f(1)=1. 
\end{aligned} 
\]
Then
\[
\widetilde{B}_n(f)(x)=5x(1-x)^5+14x^2 (1-x)^4+19x^3 (1-x)^3+14x^4 (1-x)^2+5x^5 (1-x)+x^6.
\]
Its derivative is
\[
(\widetilde{B}_n(f))'(x)=5(1-x)^5+3x(1-x)^4+x^2 (1-x)^3-x^3 (1-x)^2-3x^4 (1-x)+x^5
\]
and $(\widetilde{B}_n(f))'(7/10)=-73/2000$.

It seems that it is quite difficult to construct a monotone function $f$, for which $\widetilde{B}_n(f)$ or $\widehat{B}_n(f)$ are not monotone, by means of elementary functions.
\end{ex}

In the next example we consider the sufficient condition stated in \thmref{thm1m}.  

\begin{figure}
\centering
\begin{minipage}{.5\textwidth}
  \centering
	\includegraphics[width=\textwidth]{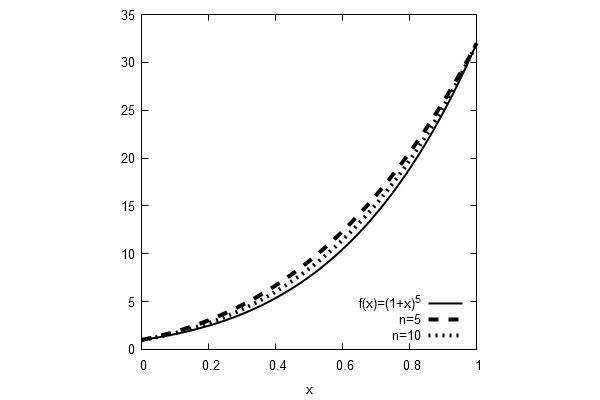}
	\caption{$\widetilde{B}_n$ and monotonicity}
	\label{fig:1}
\end{minipage}%
\begin{minipage}{.5\textwidth}
  \centering
	\includegraphics[width=\textwidth]{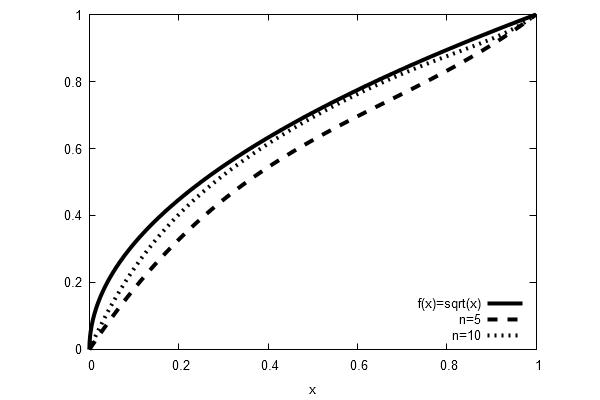}
	\caption{$\widetilde{B}_n$ and convexity}
	\label{fig:2}
\end{minipage}
\end{figure}

\begin{ex}
The function $f(x)=(x+1)^5$ satisfies the assumptions in \thmref{thm1m}. Thus the polynomials $\widetilde{B}_n(f)$ are monotone increasing for all $n$. Figure \ref{fig:1} contains the plot of $f(x)$ and $\widetilde{B}_n(f)$ for $n=5$ and $n=10$.

\end{ex}

Finally, let us demonstrate that $\widetilde{B}_n$ preserves asymptotically convexity.

\begin{ex}
Consider the concave function $f(x)=\sqrt{x}$. Figure \ref{fig:2} shows the plots of $f(x)$ and $\widetilde{B}_n(f)$ for $n=5$ and $n=10$. We can see that the graphs of $\widetilde{B}_5(f)$ and $\widetilde{B}_{10}(f)$ have an inflection point. It moves to $1$ as $n$ increases. This example shows that generally $\widetilde{B}_n$, and similarly $\widehat{B}_n$, does not preserve convexity.
\end{ex}

The computations and the plots were made with wmMaxima 16.04.2.

\bigskip
\textbf{Acknowledgements}. I am thankful to the Referee for the corrections and suggestions---they improved the presentation. Especially, I owe the Referee the elegant idea how to reduce the case of decreasing functions to the case of increasing.

\bigskip
\begin{footnotesize}
	\noindent
	\begin{tabular}{ll}
		Borislav R. Draganov& \\
		Dept. of Mathematics and Informatics&
		Inst. of Mathematics and Informatics\\
		Sofia University ``St. Kliment Ohridski''&
		Bulgarian Academy of Sciences\\
		5 James Bourchier Blvd.&
		bl. 8 Acad. G. Bonchev Str.\\
		1164 Sofia&
		1113 Sofia\\
		Bulgaria&
		Bulgaria\\
		bdraganov@fmi.uni-sofia.bg& \\
	\end{tabular}

\end{footnotesize}

\end{document}